\documentclass[a4paper,11pt]{amsart}

\usepackage[dvips]{graphics}
\usepackage[matrix,arrow,tips,curve]{xy}
\usepackage[english]{babel}
\usepackage{amsmath}
\usepackage{amssymb}

\usepackage{mathrsfs}

\usepackage{enumerate}

\usepackage{psfrag}
 \psfrag{E}{\fontsize{20}{20}$E$}
 \psfrag{F}{\fontsize{20}{20}$F$}
 \psfrag{B}{\fontsize{20}{20}$\Delta_0$}
 \psfrag{l}{\fontsize{20}{20}$\ell$}
 \psfrag{E5}{\fontsize{20}{20}$\wi{E}$}
 \psfrag{E0}{\fontsize{20}{20}$E_0$}
 \psfrag{E1}{\fontsize{20}{20}$E_1$}
 \psfrag{E2}{\fontsize{20}{20}$E_2$}
 \psfrag{A1}{\fontsize{20}{20}$A_1$}
 \psfrag{A2}{\fontsize{20}{20}$A_2$}
 \psfrag{G1}{\fontsize{20}{20}$E_0'$}
 \psfrag{G}{\fontsize{20}{20}$E_0'$}
 \psfrag{G2}{\fontsize{20}{20}$E_0''$}
 \psfrag{D1}{\fontsize{20}{20}$D_1$}
 \psfrag{D}{\fontsize{20}{20}$D$}
 \psfrag{D2}{\fontsize{20}{20}$D_2$}
 \psfrag{R1}{\fontsize{20}{20}$R_1$}
 \psfrag{R2}{\fontsize{20}{20}$\w{R}_2$}
 \psfrag{R3}{\fontsize{20}{20}$\w{R}_3$}
 \psfrag{X3}{\fontsize{20}{20}$X_3$}
 \psfrag{R4}{\fontsize{20}{20}$\w{R}_4$}
 \psfrag{E3}{\fontsize{20}{20}$E_3$}
 \psfrag{E4}{\fontsize{20}{20}$E_4$}
 \psfrag{S3}{\fontsize{20}{20}$\sigma(E_3)$}
 \psfrag{S4}{\fontsize{20}{20}$\sigma(E_4)$}

\newtheorem{thm}{Theorem}[section]

\newtheorem*{thm*}{Theorem}

\theoremstyle{definition}
\newtheorem{remark}[thm]{Remark}
\newtheorem{parg}[thm]{}

\newcommand{\ph}{\varphi}
\newcommand{\w}{\widetilde}
\newcommand{\ma}{\mathcal}
\newcommand{\la}{\longrightarrow}
\newcommand{\wi}{\widehat}
\newcommand{\pr}{\mathbb{P}}

\newcommand{\R}{\mathbb{R}}
\newcommand{\Z}{\mathbb{Z}}

\newcommand{\Gr}{\operatorname{Gr}}

\newcommand{\Pic}{\operatorname{Pic}}

\newcommand{\im}{\operatorname{Im}}
\newcommand{\NE}{\operatorname{NE}}

\newcommand{\Eff}{\operatorname{Eff}}

\newcommand{\Id}{\operatorname{Id}}

\title[Quasiparabolic vector bundles on $\pr^1$]{Rank $2$ quasiparabolic vector bundles on $\pr^1$ and the variety of linear subspaces contained in two odd-dimensional quadrics}
\author{C.~Casagrande}
\address{Universit\`a di Torino,
Dipartimento di Matematica,
via Carlo Alberto 10,
10123 Torino - Italy}
\email{cinzia.casagrande@unito.it}
\date{October 12, 2014}

\begin{document}
\maketitle

\renewcommand{\theequation}{\thethm}

\noindent Consider $2g+1$ distinct points
 $p_1,\dotsc,p_{2g+1}\in\pr^1$, where $g\geq 2$.
We recall that a rank $2$ quasiparabolic vector bundle on the marked curve 
$(\pr^1,p_1,\dotsc,p_{2g+1})$ is a rank $2$ vector bundle $V$  on $\pr^1$, with the additional data of a one-dimensional subspace $F_{j}$ of the fiber $V_{p_j}$ of $V$ over $p_j$, for every $j=1,\dotsc,2g+1$. 
The notion of stability for quasiparabolic bundles usually depends on the choice of some weights. In this paper we will only consider stability with respect to the weights $\{0,\frac{1}{2}\}$ at each marked point, and we will say that a quasiparabolic vector bundle is stable if it stable with respect to these weights (see \S\ref{mod} and references therein).
By \cite{mehtaseshadri}, there is  a fine, projective moduli space $\ma{N}$ of stable quasiparabolic vector bundles of rank $2$ and degree zero on $(\pr^1,p_1,\dotsc,p_{2g+1})$.

The purpose of this note is to show that $\ma{N}$ is isomorphic to the variety of $(g-2)$-dimensional linear subspaces of $\pr^{2g}$ contained in the intersection of two quadrics. 
\begin{thm*}
\label{main}
Let 
 $p_1,\dotsc,p_{2g+1}\in\pr^1$ be distinct points, with $g\geq 2$; assume that 
 $p_j=(\lambda_j:1)$ with $\lambda_j\in k$. Consider $\pr^{2g}$ with homogeneous coordinates $(x_1:\cdots:x_{2g+1})$, and 
let $Q_1$ and $Q_2$ denote the following quadrics:
$$Q_1\colon \sum_{j=1}^{2g+1}x_j^2=0,\qquad Q_2\colon  \sum_{j=1}^{2g+1}
\lambda_jx_j^2=0.$$
Then 
the moduli space $\mathcal{N}$ of stable quasiparabolic vector bundles of rank $2$ and degree zero on $(\pr^1,p_1,\dotsc,p_{2g+1})$
is isomorphic to the variety $\ma{G}$ of $(g-2)$-dimensional linear subspaces of $\pr^{2g}$, contained in $Q_1\cap Q_2$.
\end{thm*}
\noindent Notice that the choice  of the degree is not relevant here, as for every $d\in\Z$ $\ma{N}$ is isomorphic to the moduli space  of stable quasiparabolic vector bundles of rank $2$ and degree $d$ on $(\pr^1,p_1,\dotsc,p_{2g+1})$, see \S\ref{mod}.

The proof of this result relies on two well-known facts. The first is the relation between quasiparabolic vector bundles on $\pr^1$ and invariant vector bundles on hyperelliptic curves, established by Bhosle \cite{bhosle84}.
The second ingredient is the description by Bhosle and Ramanan \cite{bhosleramanan} of the moduli space of stable vector bundles on a hyperelliptic curve of genus $g$, with rank $2$ and fixed determinant of odd degree, as the variety of $(g-2)$-dimensional linear subspaces  of $\pr^{2g+1}$ contained in the intersections of two quadrics.

Let us notice that the variety $\ma{G}$ is a remarkable example of Fano variety; it has dimension $2g-2$, Picard number $\rho_{\ma{G}}=2g+2$, and $-K_{\ma{G}}$ very ample (see \S\ref{grass}). Both varieties $\ma{N}$ and $\ma{G}$ have been studied in several papers, see for instance \cite{bauer,biswas,abe,mukaiADE,BHK} for $\ma{N}$, and \S\ref{grass} and references therein for $\ma{G}$.

The moduli space $\ma{N}$ has a rich birational geometry: it has been shown by Bauer \cite{bauer} that it is a small modification of the blow-up of $\pr^{2g-2}$ in $2g+1$ points, see \S\ref{bir}.
In particular, we deduce that $\ma{G}$ is a rational variety.

Summing-up, we have
three different descriptions for the same Fano variety: an embedded description in a grassmannian, a modular description via quasiparabolic vector bundles on $\pr^1$, and a birational description as the unique Fano small modification of the blow-up of $\pr^{2g-2}$ in $2g+1$ points.

\medskip

\noindent{\bf Acknowledgments. } This work grew out of an attempt to construct examples of Fano $4$-folds as small modifications of blow-ups of $\pr^4$ at finitely many points. I am very grateful to Alexander Kuznetsov and Sergey Galkin for pointing out to me K\"{u}chle's work \cite{kuchle} on Fano $4$-folds, in particular the example of the Fano $4$-fold $\ma{G}$ of lines contained in the intersection of two quadrics in $\pr^6$. I would also like to thank Ana-Maria Castravet who informed me about Bauer's paper \cite{bauer} and the account of it in Mukai's book \cite[\S 12.5]{mukaibook}, 
and St\'ephane Druel for useful comments.

This research was started while I  was a guest at the Max Planck Institute for Mathematics in Bonn, in spring 2014. It is a pleasure to thank the Institute for the support and for the wonderful working conditions.
\section{Preliminaries}
\stepcounter{thm}
\subsection{Notations}
If $V$ is a vector bundle over a curve $C$, we denote by $V_p$ the fiber of $V$ over $p\in C$. Moreover, if $\alpha\colon V'\to V$ is a homomorphism between locally free sheaves, we denote by $\alpha_p\colon V'_p\to V_p$ the induced linear map.

We denote by $\text{Gr}(r,\pr^n)$ the grassmannian of $r$-dimensional linear subspaces of $\pr^n$.

If $\ma{M}$ is a moduli space, we denote by $[E]$ the point of $\ma{M}$ corresponding to the object $E$.

We work over an algebraically closed field $k$ of characteristic zero.
\stepcounter{thm}
\subsection{Quasiparabolic vector bundles}\label{exseq}
Let $C$ be a smooth projective curve with distinct marked points $p_1,\dotsc,p_r\in C$. A rank $2$ quasiparabolic vector bundle on $(C,p_1,\dotsc,p_r)$ is given by  $(V,F_1,\dotsc,F_{r})$, where $V$ is a rank $2$ vector bundle on $C$, and $F_j$ is a one-dimensional subspace of $V_{p_j}$ for every $j=1,\dotsc,r$.

Let us describe a well-known construction for shifting degrees of quasiparabolic vector bundles (see 
\cite[Rem.~5.4]{mehtaseshadri} and  \cite[\S 12.5]{mukaibook}).
 Given $(V,F_1,\dotsc,F_{r})$,
consider the natural sheaf map
$$\beta\colon V\la\bigoplus_{j=1}^{r} (V_{p_j}/F_j)\otimes\mathcal{O}_{p_j},$$ and set $V':=\ker\beta$, so that 
$V'$ is the subsheaf of sections $s$ of $V$ such that $s(p_j)\in F_j$ for all $j$.
We have an exact sequence of sheaves on $C$:
$$0\la V'\stackrel{\alpha}{\la}V\stackrel{\beta}{\la} \bigoplus_{j=1}^{r} (V_{p_j}/F_j)\otimes\mathcal{O}_{p_j} \la 0.$$
Then $V'$ is locally free of rank $2$, $\det V'\cong\det V\otimes\ma{O}_C(-p_1-\cdots-p_{r})$ (hence $\deg V'=\deg V-r$),
and $\alpha$ is an isomorphism outside $p_1,\dotsc,p_{r}$, while we have $\im\alpha_{p_j}=F_j$ and $\dim\ker\alpha_{p_j}=1$ for every $j=1,\dotsc,r$.
Thus we get a new rank $2$ vector bundle $V'$ on $C$, with a quasiparabolic structure at the points $p_1,\dotsc,p_r$ given by the one-dimensional subspaces $F_j':=\ker\alpha_{p_j}$. 

Starting from  $(V',F_1',\dotsc,F_r')$, we can repeat the procedure and get a new exact sequence
$$0\la V''\stackrel{\alpha'}{\la}V'\stackrel{\beta'}{\la} \bigoplus_{j=1}^{r} (V_{p_j}'/F_j')\otimes\mathcal{O}_{p_j} \la 0$$
and a new rank $2$ quasiparabolic vector bundle $(V'',F_1'',\dotsc,F_r'')$. Then 
$V''$ is the subsheaf of sections $s$ of $V$ vanishing at all $p_j$'s, namely
$V''\cong V\otimes\mathcal{O}_C(-p_1-\cdots-p_r)$, and the subspaces $F_j''$ correspond to $F_j$ under this isomorphism.
This shows that the quasiparabolic vector bundles  $(V,F_1,\dotsc,F_{r})$ and 
 $(V',F_1',\dotsc,F_r')$ determine each other.
\stepcounter{thm}\subsection{Moduli of stable quasiparabolic vector bundles}\label{mod}
A rank $2$ quasiparabolic vector bundle $(V,F_1,\dotsc,F_{2g+1})$ on $(\pr^1,p_1,\dotsc,p_{2g+1})$ is \emph{stable} (with respect to the weights $\{0,\frac{1}{2}\}$ at each $p_j$) if for every line subbundle $L\subset V$ we have:
$$\#\left\{j\in\{1,\dotsc,2g+1\}\,|\,L_{p_j}=F_j\right\}< \deg V+g-2\deg L+\frac{1}{2}$$
(see \cite[Def.~12.45 and Def.~12.55]{mukaibook}; in Mukai's notation, the weight at each point is $\frac{1}{2}$).

Notice that as the right-hand side is not an integer, we can never get equality above; this depends on the
fact that the number of marked points is odd, and corresponds to the non-existence of  strictly semistable quasiparabolic vector bundles in our setting.

There exists a smooth, projective, fine moduli space $\mathcal{N}_{d}$ for stable quasiparabolic vector bundles of degree $d$ and rank $2$ on $(\pr^1,p_1,\dotsc,p_{2g+1})$ \cite{mehtaseshadri}.

As usual, by tensoring with a line bundle, one sees that $\ma{N}_d\cong\ma{N}_{d'}$ when $d$ and $d'$ have the same parity. Moreover, using the construction described in \S\ref{exseq} (see also \cite[Lemma 12.51]{mukaibook}), we also get $\ma{N}_d\cong\ma{N}_{d-2g-1}$. We conclude that the moduli spaces $\ma{N}_d$ are isomorphic for all $d$; in the rest of the paper we will set $d=0$ and just write $\ma{N}$ for $\ma{N}_0$. In this case the stability condition becomes:
\stepcounter{thm}
\begin{equation}\label{stab}
\#\left\{j\in\{1,\dotsc,2g+1\}\,|\,L_{p_j}=F_j\right\}\leq g-2\deg L
\end{equation}
for every line subbundle $L\subset V$. 
\stepcounter{thm}\subsection{The variety $\ma{G}$ of $(g-2)$-dimensional linear subspaces contained in two general quadrics in $\pr^{2g}$}\label{grass}
Let $Z\subset\pr^{2g}$ be the complete intersection of two quadrics. It is well-known (see \cite[Prop.~2.1]{reidthesis} or \cite[Lemma 8.6.1]{dolgachevbook}) that $Z$ is smooth if and only if, up to a projective transformation of $\pr^{2g}$, we have $Z=Q_1\cap Q_2$
where $Q_1$ and $Q_2$ are the quadrics: $$Q_1\colon \sum_{j=1}^{2g+1}x_j^2=0\qquad\text{and}\qquad Q_2\colon\sum_{j=1}^{2g+1}\lambda_jx_j^2=0,$$
with $\lambda_j\in k$ all distinct. 

Let us assume that $Z$ is indeed smooth, and let $\ma{G}$ be the variety of $(g-2)$-dimensional linear subspaces contained in $Z$.
Then $\ma{G}$ is smooth, connected, and has dimension $2g-2$ (see \cite[Th.~2.6]{reidthesis} for smoothness, and \cite[Th.~4.1]{borcea90} or \cite[Th.~2.1]{debarremanivel} for connectedness). Moreover $-K_{\ma{G}}$ is given by the restriction to $\mathcal{G}$ of $\ma{O}(1)$ on the grassmannian $\text{Gr}(g-2,\pr^{2g})$ (see \cite[Rem.~4.3]{borcea90} or \cite[Rem.~3.2(2)]{debarremanivel}), in particular $\ma{G}$ is a Fano variety,  $-K_{\ma{G}}$ is very ample, and the Pl\"ucker embedding of $\text{Gr}(g-2,\pr^{2g})$ in $\pr^{\binom{2g+1}{g-1}-1}$ yields an anticanonical embedding for $\ma{G}$.
Finally we have $\rho_{\ma{G}}=2g+2=\dim\ma{G}+4$ \cite[Prop.~3.2]{jiang}.
\section{Proof of the Theorem}
\begin{parg}
As a first step, we embed $\ma{G}$ in the variety $\ma{G}'$ of $(g-2)$-dimensional linear subspaces of $\pr^{2g+1}$, contained in the intersection of two $(2g)$-dimensional quadrics.

Consider $\pr^{2g+1}$ with homogeneous coordinates $(x_1:\cdots:x_{2g+2})$. We identify $\pr^{2g}$ with the hyperplane $H:=\{x_{2g+2}=0\}\subset\pr^{2g+1}$, and $\Gr(g-2,\pr^{2g})$ with the subvariety 
$$\left\{[L]\in\Gr(g-2,\pr^{2g+1})\,|\,L\subset H\right\}.$$

Fix $\lambda_{2g+2}\in k$ different from $\lambda_1,\dotsc,\lambda_{2g+1}$,
and consider the two quadrics in $\pr^{2g+1}$: 
$$Q'_1\colon \sum_{j=1}^{2g+2}x_j^2=0\qquad\text{and}\qquad Q'_2\colon  \sum_{j=1}^{2g+2}
\lambda_jx_j^2=0.$$
Finally let $\mathcal{G}'\subset\Gr(g-2,\pr^{2g+1})$ be the variety of  $(g-2)$-dimensional linear subspaces of 
 $\pr^{2g+1}$ contained in $Q_1'\cap Q_2'$; we have 
$$\mathcal{G}=\{[L]\in \mathcal{G}'\,|\,L\subset H\}.$$

Let us consider the involutions  of  $\pr^{2g+1}$ and of $\ma{G}'$ given by:  
\stepcounter{thm}
\begin{gather}
i_{\pr^{2g+1}}(x_1\cdots:x_{2g+2})=
(x_1:\cdots:x_{2g+1}:-x_{2g+2}),\notag\\ i_{\ma{G}'}([L])=[i_{\pr^{2g+1}}(L)].\label{involution}
\end{gather}
If $L$ is a linear subspace of $\pr^{2g+1}$, an elementary computation shows that $i_{\pr^{2g+1}}(L)=L$ if and only if either $L\subseteq H=\{x_{2g+2}=0\}$, or $L$ contains the point $(0:\cdots:0:1)$. Since this point is not contained in the quadric $Q'_1$, we deduce that $\ma{G}$ is the fixed locus of the involution $i_{\ma{G}'}$.
\end{parg}
\begin{parg}
Set $p_{2g+2}:=(\lambda_{2g+2}:1)\in\pr^1$, and notice that the points $p_1,\dotsc,p_{2g+2}$ are distinct. Let $\pi\colon X\to \pr^1$ be the double cover of $\pr^1$ ramified over $p_1,\dotsc,p_{2g+2}$, so that $X$ is a hyperelliptic curve of genus $g$. Set $w_j:=\pi^{-1}(p_j)$ for $j=1,\dotsc,2g+2$, and let $i\colon X\to X$ be the hyperelliptic involution.

Let  $\mathcal{M}$ be the moduli space of stable rank $2$ vector bundles on $X$, with determinant  $\mathcal{O}_X(-w_1-\cdots-w_{2g+1})$;
by \cite{bhosleramanan} there exists an isomorphism 
$$\ph\colon \mathcal{M}\la \mathcal{G}'.$$
\end{parg}
\begin{parg}
A vector bundle $E$ on $X$ is \emph{$i$-invariant} if $i^*E\cong E$.
As $i^*\mathcal{O}_X(-w_1-\cdots-w_{2g+1})\cong \mathcal{O}_X(-w_1-\cdots-w_{2g+1})$, $i$ induces an involution $i_{\ma{M}}$ of $\ma{M}$ by sending $[E]$ to $[i^*E]$. We denote by $\mathcal{M}^{\text{inv}}$ the locus of $i$-invariant vector bundles in $\ma{M}$, namely the fixed locus of $i_{\ma{M}}$.

Set $D:=w_1+\cdots+w_{2g+1}$ and $\beta:=\mathcal{O}_X(D)$. 
 In the notation of
  \cite[p.~161]{bhosleramanan}
we have $i_{\ma{M}}=i_{\beta}$, where $i_{\beta}$ is the involution of $\ma{M}$ defined by $i_{\beta}([E])=[i^*E\otimes \mathcal{O}_X(-w_1-\cdots-w_{2g+1})\otimes\beta]$.

As in  \cite[Lemma 2.1]{bhosleramanan}, the line bundle $\beta$ corresponds to a partition $\{1,\dotsc,2g+2\}=S\cup T$,
where 
\begin{align*}
S:= &\{j\,|\,\text{the coefficient of $w_j$ in $D$ is odd}\ \}=\{1,\dotsc,2g+1\},
\\
T:=&\{j\,|\,\text{the coefficient of $w_j$ in $D$  is even}\}=\{2g\}.
\end{align*}
Notice that by choosing another divisor $D'$ linearly equivalent to $D$, and with support contained in $\{w_1,\dotsc,w_{2g+2}\}$, we get the same partition, with at most $S$ and $T$ interchanged.

By \cite[Corollary, p.~161]{bhosleramanan} $i_{\ma{M}}$ corresponds, under the isomorphism $\ph$, to the involution of $\ma{G}'$ induced by the involution of $\pr^{2g+1}$ which changes the sign of the coordinates $x_j$ for  $j\in S$. This is precisely the involution $i_{\ma{G}'}$ in \eqref{involution}. 

 We conclude that $\ph$ restricts to an isomorphism between $\mathcal{M}^{\text{inv}}$ and $\ma{G}$; in particular, $\ma{M}^{\text{inv}}$ is smooth, irreducible, and has dimension $2g-2$ (see \S\ref{grass}). We are left to show that $\mathcal{M}^{\text{inv}}$ is isomorphic to $\mathcal{N}$.  
\end{parg}
\begin{parg}\label{map}
The isomorphism $\mathcal{M}^{\text{inv}}\cong\mathcal{N}$ follows basically from \cite[Prop.~1.2]{bhosle84} (see also \cite[Prop.~3.2]{bhosle90}); we report the details for the reader's convenience.

We first describe a set-theoretical map from $\mathcal{N}$ to $\mathcal{M}^{\text{inv}}$.

Let $(V,F_1,\dotsc,F_{2g+1})$ be a rank $2$ stable quasiparabolic vector bundle on $(\pr^1,p_1,\dotsc,p_{2g+1})$, of degree zero. Its pull-back $\pi^*V$ 
inherits the quasiparabolic structure $(\pi^*V,\pi^*F_1,\dotsc,\pi^*F_{2g+1})$
on the curve $X$ with marked points $w_1,\dotsc,w_{2g+1}$.
As described in \S\ref{exseq}, we have an exact sequence of sheaves on $X$:
$$0\la E\stackrel{\alpha}{\la}\pi^*V\stackrel{\beta}{\la} \bigoplus_{j=1}^{2g+1} (\pi^*V_{w_j}/\pi^*F_j)\otimes\mathcal{O}_{w_j} \la 0,$$
where $E$ is locally free of rank $2$, $\det E\cong\ma{O}_X(-w_1-\cdots-w_{2g+1})$,
and $\alpha$ is an isomorphism outside $w_1,\dotsc,w_{2g+1}$, while $\im\alpha_{w_j}=\pi^*F_j$ and $\dim\ker\alpha_{w_j}=1$ for every $j=1,\dotsc,2g+1$.

The natural isomorphism $\pi^*V\cong i^*\pi^*V$ and the exact sequence above induce a natural isomorphism 
$$\xi\colon E \la i^*E$$ such that $i^*(\xi)\circ \xi=\Id_E$ ($\xi$ can be thought as a lifting of the involution $i$ to the total space of the vector bundle $E$); in particular, for every $j=1,\dotsc,2g+2$ we have an induced involution $\xi_{w_j}\colon E_{w_j}\to E_{w_j}$.

 As the homomorphism $\beta$ is trivial in $w_{2g+2}$, $\xi_{w_{2g+2}}$ is the identity. A local computation shows that for $j=1,\dotsc,2g+1$, $\xi_{w_j}$ has eigenvalues $1$ and $-1$, and the $(-1)$-eigenspace is $\ker\alpha_{w_j}$.
\end{parg}
\begin{parg}
Let us show that $E$ is stable (see also
 \cite[Prop.~1.2]{bhosle84}). 
By contradiction, suppose the contrary. Then there exists a line subbundle $L\subset E$ such that 
$$\deg L\geq \frac{\deg E}{2}=-g-\frac{1}{2}.$$ 
Since a rank $2$ vector bundle has at most one destabilising line bundle (see \cite[Prop.~10.38]{mukaibook}), we must have $\xi(L)=i^*L\subset i^*E$.

Let $M\subset\pi^*V$ be the line subbundle generated by the image of $L$ under $\alpha\colon E\to\pi^*V$. Since $\xi(L)=i^*L$, we must have $\wi{\xi}(M)=i^*M\subset i^*\pi^*V$ under the natural isomorphism
$ \wi{\xi}\colon \pi^*V\to i^*\pi^*V$.
 This implies that $M=\pi^*(M')$, where $M'$ is a line subbundle of $V$.

Set $J:=\{j\in\{1,\dotsc,2g+1\}\,|\,L_{w_j}=\ker\alpha_{w_j}\}$ and $J^c:=\{1,\dotsc,2g+1\}\smallsetminus J$.
 Then $L\cong M\otimes\ma{O}_X(-\sum_{j\in J}w_j)$, thus
$\deg L=\deg M-|J|=2\deg M'-|J|$, which yields
$$|J^c|=2g+1-|J|=2g+1+\deg L-2\deg M'\geq g+\frac{1}{2}-2\deg M'.$$
Notice that if $j\in J^c$, then $M_{w_j}=\im\alpha_{w_j}$, hence $M'_{p_j}=F_j$. Thus the equation above contradicts the stability of $(V,F_1,\dotsc,F_{2g+1})$ (see \eqref{stab}). 

Therefore $E$ is stable, and $[E]\in\ma{M}^{\text{inv}}$.
\end{parg}
\begin{parg}
The construction in \ref{map} can be made in families, starting from the universal family over $\ma{N}$; this yields a morphism 
$$\psi\colon\ma{N}\la\ma{M}^{\text{inv}}.$$
As $\ma{N}$ and $\ma{M}^{\text{inv}}$ are smooth, irreducible varieties of the same dimension, to conclude that  $\psi$ is an isomorphism it is enough to show
that  $\psi$ is injective.

Let  $[(V,F_1,\dotsc,F_{2g+1})]$ and $[(\w{V},\w{F}_1,\dotsc,\w{F}_{2g+1})]$ be two points of  $\ma{N}$,
with the same image $[E]\in \ma{M}^{\text{inv}}$ under $\psi$.
By construction we have two isomorphisms
$\xi,\w{\xi}\colon E\to i^*E$ such that  
\begin{gather}
\stepcounter{thm}
i^*(\xi)\circ \xi=\Id_E,\qquad
i^*(\w{\xi})\circ \w{\xi}=\Id_E,\label{first}\\ \stepcounter{thm}
\text{ and }\qquad \xi_{w_{2g+2}}=\w{\xi}_{w_{2g+2}}=\Id_{E_{w_{2g+2}}}. \label{second}
\end{gather}
As $E$ is stable, it has only constant automorphisms, and there exists a non-zero constant $\lambda$ such that $\w{\xi}=\lambda\xi$. Then \eqref{first} implies $\lambda=\pm 1$, and \eqref{second} yields $\lambda=1$, namely $\w{\xi}=\xi$.

In particular, the $(-1)$-eigenspaces of 
 $\xi_{w_j}$ and $\w{\xi}_{w_j}$ are the same for all $j=1,\dotsc,2g+1$, thus the quasiparabolic vector bundles $$(E,\ker\alpha_{w_1},\dotsc,\ker\alpha_{w_{2g+1}})\quad\text{and}\quad
(E,\ker\w{\alpha}_{w_1},\dotsc,\ker\w{\alpha}_{w_{2g+1}})$$ coincide. 
As noticed 
 in \S\ref{exseq}, this shows that the quasiparabolic vector bundles
$(\pi^*V,\pi^*F_1,\dotsc,\pi^*F_{2g+1})$ and
$(\pi^*\w{V},\pi^*\w{F}_1,\dotsc,\pi^*\w{F}_{2g+1})$
are isomorphic, and hence 
the same holds for $(V,F_1,\dotsc,F_{2g+1})$ and $(\w{V},\w{F}_1,\dotsc,\w{F}_{2g+1})$.

 This shows that $\psi$ is injective, and concludes the proof of the Theorem.
\qed
\end{parg}
\begin{remark}
The Verlinde formula \cite[\S 12.5, in particular Remark 12.54]{mukaibook} gives 
$$h^0(\ma{G},\ma{O}(-K))=h^0(\ma{N},\ma{O}(-K))=1+4+4^2+\cdots+4^{g-1}=\frac{4^g-1}{3}.$$
On the other hand, we have
 $\mathcal{G}\subset\pr^{\binom{2g+1}{g-1}-1}$ under
the Pl\"{u}cker embedding. Then one can check that $h^0(\ma{G},\ma{O}_{\ma{G}}(-K_{\ma{G}}))=\binom{2g+1}{g-1}$ for $g=2,3$, while 
 $h^0(\ma{G},\ma{O}_{\ma{G}}(-K_{\ma{G}}))>\binom{2g+1}{g-1}$ for $g\geq 4$, so that $\ma{G}\subset \pr^{\binom{2g+1}{g-1}-1}$ is not linearly normal for $g\geq 4$; see \cite[Th.~1 and Th.~3]{kuchle96} and \cite[Rem.~4.2(2)]{debarremanivel} for related results.
Notice that instead, for $g=2$, $\ma{G}\subset\pr^4$ is a Del Pezzo surface of degree $4$, and is projectively normal (see for instance 
\cite[Th.~8.3.4]{dolgachevbook}).
\end{remark}
\begin{remark}[The case of an even number of marked points]
Let $\ma{N}^+_d$ be the moduli space of semistable 
rank two quasiparabolic vector bundles of degree $d$ on $(\pr^1,p_1,\dotsc,p_{2g+2})$. As in \S\ref{map}, one can associate to $(V,F_1,\dotsc,F_{2g+2})$ an $i$-invariant vector bundle $E$ on the hyperelliptic curve $X$; however, the degree of $E$ is even, so this relates $\ma{N}^+_d$ to the moduli space $\ma{M}^+$ of semistable rank two vector bundles on $X$ with fixed determinant of even degree. Moreover, the resulting map $\ma{N}^+_d\to\ma{M}^+$ is not injective, but has degree two onto its image; for more details we refer the reader to \cite[Th.~2.1]{kumar} and \cite[\S 2.13]{abe} and references therein. 
\end{remark}
\section{Final remarks}
\stepcounter{thm}
\subsection{Birational geometry: relation with the blow-up of $\pr^n$ at $n+3$ general points}\label{bir}
Set $n:=2g-2$, so that $2g+1=n+3$. Let $q_1,\dotsc,q_{2g+1}\in\pr^n$ be the images of $p_1,\dotsc,p_{2g+1}\in\pr^1$ under the Veronese embedding $\pr^1\hookrightarrow\pr^n$; notice that the points  $q_1,\dotsc,q_{n+3}$ are in general position, in the sense that any $n+1$ points are projectively independent.

Let $Y$ be the blow-up of $\pr^n$ at $q_1,\dotsc,q_{n+3}$. It has been shown by Bauer \cite{bauer} (see also \cite[Th.~12.56]{mukaibook}) that there exists a birational map $\ma{N}\dasharrow Y$, which is an isomorphism in codimension one. More precisely, $Y$ is a Mori dream space, $\ma{N}$ is its (unique) Fano small modification, and the birational map $\ma{N}\dasharrow Y$ factors a sequence of $K$-negative flips. 

The cone  of effective divisors  $\Eff(Y)\subset\Pic(Y)\otimes\R$ and the Cox ring of $Y$ are described in \cite{mukaiXIV,CoxCT}; $\Eff(Y)$ has $2^{n+2}$ extremal rays.  
Since the Picard group, the cone of effective divisors, and the Cox ring are invariant under the small modification  $\ma{N}\dasharrow Y$, the same description applies as well to $\ma{N}$ and $\ma{G}$.

Notice that when $g=2$, $\ma{G}$ is just the intersection of two quadrics in $\pr^4$ (namely a Del Pezzo surface of degree $4$), $Y$ is the blow-up of $\pr^2$ at $5$ points, and $\ma{G}\cong Y$.
\stepcounter{thm}
\subsection{Dimension $4$}
Let us set $g=4$ in this paragraph, so that $\ma{G}\subset\text{Gr}(1,\pr^6)$ is the variety of lines contained in the intersection of two quadrics in $\pr^6$, and
has dimension $n=4$.
The Fano $4$-fold $\ma{G}$ has been studied by Borcea \cite{borcea91} and K\"{u}chle \cite{kuchle}; it has $b_2=8$, $b_3=0$, $b_4=h^{2,2}=30$, $(-K)^4=80$, and $h^0(-K)=21$. It is a peculiar example of Fano $4$-fold, because it has ``large'' second Betti number: the only other examples known to the author of Fano $4$-folds with $b_2\geq 8$ are products of Del Pezzo surfaces.

As above let $Y$ be the blow-up of $\pr^4$ in $7$ points.
The small modification 
$Y\dasharrow \ma{G}$ has a simple, explicit description as a sequence of $22$ $K$-positive flips, see \cite[Ex.~12.57]{mukaibook}.

The cone of effective curves $\NE(\ma{G})$ has $64$ extremal rays; these are all small, of type $(2,0)$, with exceptional locus a surface $L\cong\pr^2$  with normal bundle $\ma{O}(-1)^{\oplus 2}$ \cite[Th.~4.3]{borcea91}. These surfaces in fact are given by the lines contained in the $64$ planes contained in $Q_1\cap Q_2\subset\pr^6$. 


We conclude by noting that the blow-up of $\ma{G}$ at a general point is still a Mori dream space, because it is a small modification of a blow-up of $\pr^4$ at $8$ general points, which is a Mori dream space (see \cite[Th.~1.3]{CoxCT} and also \cite[\S 2]{mukaiADE}). It would be interesting to know whether this blow-up still has a Fano small modification; this would give an example of Fano $4$-fold with $b_2=9$, which is not a product of surfaces.
\providecommand{\bysame}{\leavevmode\hbox to3em{\hrulefill}\thinspace}
\providecommand{\MR}{\relax\ifhmode\unskip\space\fi MR }
\providecommand{\MRhref}[2]{%
  \href{http://www.ams.org/mathscinet-getitem?mr=#1}{#2}
}
\providecommand{\href}[2]{#2}

\end{document}